\documentclass{article}

\usepackage[a4paper]{geometry}
\usepackage[british]{babel}
\usepackage[utf8]{inputenc}
\usepackage{amsmath,amsfonts,amssymb,euscript,enumitem}
\usepackage{hyperref}

\title{Translation of ``Simplizialzerlegungen von Beschrankter Flachheit'' by Hans Freudenthal, Annals of Mathematics, Second Series, Volume 43, Number 3, July 1942, Pages 580-583.\footnote{Freudenthal added the following as a footnote to the title: This note is essentially the same as a note that was submitted to Fundamenta Mathematicae in March 1939. However this note never appeared in print. The result of was also used in a different paper ``Die Triangulation der differenzierbaren Mannigfaltigkeiten'' Nederlandse Akademie voor Wetenschappen, Procedings, 42 (1939), 880–901, \url{https://dwc.knaw.nl/DL/publications/PU00014650.pdf}. [More extensive reference inserted by translator]
}}
\author{Translated by:\\
Mathijs Wintraecken, Inria Sophia-Antipolis, Universit{\'e} C{\^o}te d'Azur \\
and Institute of Science and Technology Austria \footnote{The translator has been supported by the European Research Council (ERC) under  the  European  Union's  Seventh  Framework  Programme (FP/2007-2013)  /  ERC  Grant Agreement No. 339025 GUDHI (Algorithmic Foundations of Geometry Understanding in Higher Dimensions), the European Union's Horizon 2020 research and innovation programme under the Marie Sk{\l}odowska-Curie grant agreement No. 754411, and the Austrian science fund (FWF) under grant agreement M-3073. }}

\date{}

\begin{document}

\maketitle

We answer a question by Brouwer about the construction of an infinite series of subdivisions of a polytope, such that the next element in the sequence is a subdivision of the previous one and such that the subsimplices that arise do not become arbitrarily flat, that is the quotient 
\[c^r/ v\] 
($c=$ diameter $=$ the longest edge of the simplex,  $v =$ volume of the simplex) is uniformly bounded for all $r$ dimensional simplices.
Such tillings are useful in analysis and the interface of analysis and topology.

The construction, that we give here, is analogous to our subdivision by simplices of the Cartesian product of two simplices.\footnote{Freudenthal, ``Eine Simplizialzerlegung des Cartesischen Produktes zweier Simplexe'' Fundamenta Mathematicae, 29 (1937), 138-144, \url{https://bibliotekanauki.pl/articles/1383987}.}

\addcontentsline{toc}{section}{1}
\begin{center} \textbf{1} 
\end{center}

Simplices will be denoted by a fixed order of vertices, parallelepiped will be given by a marked vertex and an ordered set of edges that emanate from this vertex.

The simplex $T$ has vertices
\[ e_0 ,\dots , e_r.\]
We can also describe $T$ by the vectors
\begin{align} 
x_0 &=e_0,  & x_i &= e_i - e_{i-1} & & (i= 1 ,\dots ,r) .
\nonumber
\end{align}
Let $\pi$ be a permutation $p_1, \ldots , p_r$ of the numbers $1 \ldots r$. The points\footnote{Note translator: There was a typo in the sum. }
\begin{align} 
e_i &= x_0 + \sum_{\nu=1}^{i} x_{p_\nu}   &  &(i= 1 ,\dots ,r) 
\nonumber
\end{align} 
(seen as vertices) yield a simplex $T^\pi$. If $\pi_0$ the identity permutation is of $1 , \dots, r$ then $T^{\pi_0} = T$. The $r!$ simplices $T^\pi$ are said to the conjugates of each other.

$T^\pi$ is the set\footnote{Note translator: Here I updated the notation quite significantly. }
\[T^\pi =  \left \{  \sum_{\nu=1}^{r} \lambda_ \nu e_{ p_\nu}  \middle| \lambda_\nu \geq 0 , \sum_{\nu =1}^r \lambda_\nu =1  \right\}, \] 
or what boils down to the same thing,\footnote{Note translator: This is not entirely straightforward, the argument (or a slight generalization thereof) is spelled out in \cite{boissonnat_et_al:LIPIcs.SoCG.2021.17}. } 
\[ T^\pi= \left \{  x_0 + \sum_{\nu=1}^{r} \alpha_\nu x_{ p_\nu}  \middle| 1\geq \alpha_1 \geq  \dots \geq \alpha_r \geq 0  \right\}. \]
It follows that: The $T^\pi$ are the images of a subdivision of the parallelpiped $P$, given by, 
\[ P=\left \{  x_0 + \sum_{\nu=1}^{r} \alpha_\nu x_{ \nu}  \middle| 1\geq \alpha_\nu \geq 0  \right\}. \]

Conversely the simplex $T^\pi$ is unambiguously determined by the parallelepiped $P$, that is given by the vertex $x_0$ and the edge vectors $x_1, \dots, x_r$ and by the permutation $\pi$.

\addcontentsline{toc}{section}{2}
\begin{center} \textbf{2} \end{center} 

Cutting all the edges of the parallelepiped in half all, yields [note translator: $2^r$] parallelepiped 
\[P_\sigma =\left \{  x_{\sigma, 0} + \frac{1}{2} \sum_{\nu=1}^{r} \alpha_\nu x_{ \nu}  \middle| 1\geq \alpha_\nu \geq 0  \right\}  \] 
where\footnote{Note translator: Here I updated the notation quite significantly. } $\sigma= (s_1 ,\dots ,s_r)  \in \{ 0,1\}, ^r$, 
and
\[ 
x_{\sigma,0}= x_0 +\frac{1}{2} \sum_{\nu =1}^r s_\nu x_\nu 
\]
In the same way that $P$ is subvidived by $T^\pi$, is $P_\sigma$ subdivided by $(T_\sigma)^\pi$. Note that $T^\pi$ and $(T_\sigma )^\pi$ are homothetic [with the ratio] $(1:2)$.

Let $\pi$ again be a permutation $p_1\dots  p_r$ of $1 , \dots,r$ and let $\pi^* = \pi ^\sigma$ be the permutation $p_1^* \dots p_r^* $ that arises if one first (using the order given by $\pi$) takes all $p_\nu$ with $s_{p_\nu}=1$ and then all with $s_{P_\nu}$, so
\begin{itemize}
\item if $s_a> s_b$ then $a$ comes before $b$ in $\pi^*$
\item if $s_a = s_b$ and $a$ is before $b$ in $\pi$ then $a$ also comes before $b$ in $\pi^*$. 
\end{itemize} 
We claim that $(T_\sigma )^\pi$ lies completely in $T^{\pi^\sigma}$.
Indeed we have, 
\begin{align} 
(T_\sigma )^\pi &=  \left \{  x_0 + \frac{1}{2} \sum_{\nu=1}^{r} s_\nu x_\nu + \frac{1}{2} \sum_{\nu=1}^{r} \alpha_\nu x_{p_\nu}  \middle| 1\geq \alpha_1 \geq  \dots \geq \alpha_r \geq 0   \right\}
\nonumber 
\\
&= 
\left \{  x_0 + \sum_{\nu\in \{ \nu | s_{p_\nu} =1 \} } \frac{1}{2} (\alpha_\nu +1) x_{p_ \nu}  + \frac{1}{2} \sum_{\nu\in \{ \nu | s_{p_\nu} =0 \} } \alpha_\nu x_{ \nu}  \middle| 1\geq \alpha_1 \geq  \dots \geq \alpha_r \geq 0   \right\}
\nonumber
\\
&= \left \{ x_0 + \sum_{\nu =1}^{r} \beta_\nu x_{p_\nu^{*}} \middle| 1 \geq \beta_1 \geq \dots \geq \beta_u \geq \frac{1}{2} \geq \beta_{u+1} \geq \dots \geq \beta_r \geq 0 \right\} ,
\nonumber
\end{align} 
where $u$ is $\# \{ \nu | s_{p_\nu} =1 \}$. Our claim now follows.
Moreover, we have that the $(T_\sigma)^\pi$ such that $\pi ^\sigma = \pi^*$  form a simplicial subdivision $Z(T)$ of $T$.

\addcontentsline{toc}{section}{3}
\begin{center} \textbf{3}\footnote{Note translator: I personally find the subdivision easier to see, by first observing that it is a hyperplane arrangement and then do a rescaling, compare to \cite{boissonnat_et_al:LIPIcs.SoCG.2021.17}, also see the references mentioned there.} \end{center}

Let $e_{ij}^{\pi^*}$ be the midpoint of the edge $e_i^{\pi^*} e_j^{\pi^*}$, where we define  $e_{ii}^{\pi^*}=e_i^{\pi^*}$. With this definition, and assuming without loss of generality that $i \leq j$, we have\footnote{Note translator: Here too there was a typo, I doubt that Freudenthal ever saw the proofs of this paper, because of the  Nazi occupation of the Netherlands and he being Jewish.}
\begin{align}
 e_{ij}^{\pi^*} &= \frac{1}{2} \left ( x_0 + \sum_{\nu =1}^{i}  x_{p_\nu^*}  + x_0 + \sum_{\nu =1}^{j}  x_{p_\nu^*} \right) 
\nonumber 
\\ 
&= x_0 +\sum_{\nu =1}^{i}  x_{p_\nu^*}  + \frac{1}{2} \sum_{\nu =i+1}^{j}  x_{p_\nu^*} 
\nonumber
\end{align} 
We define $k'$ to be the number of $\nu$ with $s_{p_\nu}=1$, and $k''$ to be the number of $\nu$ with $s_{p_\nu}=0$. 

We say that the $k$th vertex of $T_\sigma ^\pi$ is given by setting $\alpha_1 =  \dots= \alpha_k =1$ and $\alpha_{k+1} = \dots = \alpha_r=0$. 
The $k$th vertex of $T_\sigma ^\pi$ can we written in terms of the midpoints we defined above as \footnote{[Literally:] this therefore gives}
\begin{align}
x_0 + \sum_{\nu=1} ^{k'} x_{p_\nu^*} + \sum_{\nu=k'+ 1}^{k'' +u} x_{p^*_\nu} = e_{ij} ^{\pi^*}, 
\nonumber
\end{align} 
where $i= k'$ and $j= k'' +u$.

From this it follows that for every subsimplex $T_\sigma ^\pi$ of $(T ^\pi)^\sigma$:
\begin{itemize} 
\item The 0th vertex is $e^{\pi ^ \sigma} _{ 0 u}$
\item the vertex $e^{\pi ^ \sigma} _{ ij}$ is followed by either $e^{\pi ^ \sigma} _{ i+1 \,j}$  or $e^{\pi ^ \sigma} _{ i\, j+1}$. 
\item the first index is bounded by u and the second by r.
\end{itemize} 
Conversely are the subsimplices $(T_\sigma)^\pi$ of $T^{\pi^*}$ characterized by these properties.

\begin{center} \textbf{4} \end{center} 

We now concentrate on $T= T^\pi_0$ and drop the index $\pi_0$ in the following. We can now also describe the subdivision $Z(T)$ of $T$ in the following way:
\begin{itemize} 
\item The vertices of $Z(T)$ are the $e_{ij}$;
\item The $e_{ij}$ and $e_{i'j'}$ form a one dimensional simplex if and only if the pairs $ij$ and $i' j'$ do not separate each other;
\item a set of $e_{ij} e_{ij}$  form a simplex of $Z(T)$ is and only if the elements form pairwise simplices of $Z(T)$;
\item the vertices in the simplices of Z(T)are ordered in ascending order of $i+j$.
\end{itemize} 
We also note that:
The simplices of $Z(T)$ are conjugated with $T$. 

Let now a finite polytope $R$ be given. We impose an ordering on the vertices and form the division $Z(R)$, where we use the process $Z$ on every simplex of $R$. The edges of $Z(R)$ are ordered lexicographically.  If we repeat the process as many times as we want we preserve a order of subdivision. All simplices that occur are similar to conjugated simplices of $R$ and because similar simplices have the same quality
\[c^r/ v\] 
we see that the flatness of the simplices remains bounded.

\textbf{Amsterdam}  

\paragraph{Acknowledgements}
The translator thanks the editorial board of the Annals of Mathematics for their permission to make a translation of this paper public. The Annals of Mathematics has the copyright to the original German text. 

\phantomsection
\bibliographystyle{alpha}
\addcontentsline{toc}{section}{Bibliography}
\bibliography{biblioShort}

\end{document}